\documentclass[12pt]{amsart}

\usepackage[utf8]{inputenc}
\usepackage{amsmath}
\usepackage{amssymb}
\usepackage{amsthm}
\usepackage{graphicx}

\newtheorem{theorem}{Theorem}

\newtheorem{corollary}{Corollary}

\title[An application of Brouwer's fixed-point theorem...]{An application of Brouwer's fixed-point theorem: continuously differentiable convex functions with gradient of constant norm}
\author{Csaba Vincze}
\address{Institute of Mathematics, University of Debrecen, H-4002 Debrecen, P. O. Box 400, Hungary}
\email{csvincze@science.unideb.hu}
\keywords{Convex functions. Gradient. Brouwer's fixed-point theorem. Browder-Minty theorem.}
\subjclass{26B25}

\begin{document}
\begin{abstract}
As an application of Brouwer's fixed-point theorem we prove that a continuously differentiable convex function with gradient of constant norm is an affine mapping. It is a first-order characterization of affine mappings among continuously differentiable convex functions, because neither the second-order condition of convexity nor related operators are used. The condition of differentiability is essential as the case of the norm function shows. In addition to Brouwer's theorem, the proof is based on the Cauchy--Bunyakovsky--Schwarz inequality and becomes complete by minimizing the distance between lines of gradient directions. 

Following the steps of the proof, we sketch a possible generalization of the result to functions defined on Hilbert spaces.
\end{abstract}
\maketitle

\section{Introduction}

Let $f\colon \mathbb{R}^n\to \mathbb{R}$ be a continuously differentiable convex function. Its gradient will be denoted by $\nabla f$. We are going to prove that if the gradient of the function is of constant norm then the function is an affine mapping of the form 
$$f(u)=\langle c_1, u \rangle+c_0 \quad (u\in \mathbb{R}^n).$$
Dropping the differentiability condition, the norm functions of the form
$$f(u)=|u|+c_0 \quad (t\in \mathbb{R}^n)$$
are obvious examples for functions with gradients of constant norm. The so-called distance functions form a more general family of singular examples \cite{Borwein} (Lem. 15.): {\emph{given a non-empty closed set $K\subset \mathbb{R}^n$, if the
distance function generated by $K$ is Fr\'{e}chet differentiable at $u\notin K$ then its gradient at $u$ is of unit norm}}. As a theoretical e\-xample, let $n\geq 2$ be a given integer and consider a continuously differentiable function $g\colon \mathbb{R}^{n-1}\to \mathbb{R}$. Taking the function 
$$f\colon \mathbb{R}^n\to \mathbb{R}, \quad f(u):=\inf_{x\in \mathbb{R}^{n-1}} \left | (x, g(x))-u\right |,$$
let us define the projection $\pi\colon \mathbb{R}^n\to \mathbb{R}^{n-1}$ such that
$$\pi(u)\in \{x^*\in \mathbb{R}^{n-1}\ | \ f(u)=\left | (x^*, g(x^*))-u\right |\}$$
for any $u\in \mathbb{R}^n$. Since $\pi(u)$ is a minimizer of the expression
$$\left | (x, g(x))-u\right |^2=\left(g(x)-u_n\right)^2+\sum_{k=1}^{n-1}(x_k-u_k)^2 \quad (x\in \mathbb{R}^{n-1}),$$
the vanishing of the partial derivatives at $x=\pi(u)$ implies that
\begin{equation}
\label{extremal}
(g\circ \pi(u)-u_n)D_k g (\pi(u))+\pi_k(u)-u_k=0 \quad (k=1, \ldots, n-1).
\end{equation}
Using that 
$$f(u)=\sqrt{\left(g\circ \pi(u)-u_n\right)^2+\sum_{k=1}^{n-1}(\pi_k(u)-u_k)^2},$$
we have the partial derivatives
$$D_i f(u)=\frac{(g\circ \pi (u)-u_n)\sum_{k=1}^{n-1}D_k g (\pi(u))D_i\pi_k(u)+\sum_{k=1}^{n-1}(\pi_k(u)-u_k)(D_i \pi_k(u)-\delta_{ik})}{f(u)}=$$
$$\frac{\sum_{k=1}^{n-1}\left((g\circ \pi (u)-u_n)D_k g (\pi(u))+\pi_k(u)-u_k \right)D_i \pi_k(u)-\delta_{ik}(\pi_k(u)-u_k)}{f(u)}\stackrel{\eqref{extremal}}{=}$$
$$-\frac{\pi_i(u)-u_i}{f(u)} \quad (i=1, \ldots, n-1),$$
$$D_n f(u)=\frac{(g\circ \pi (u)-u_n)\left(-1+\sum_{k=1}^{n-1}D_k g (\pi(u)) D_n\pi_k(u)\right)+\sum_{k=1}^{n-1}(\pi_k(u)-u_k)D_n \pi_k(u)}{f(u)}=$$
$$\frac{-(g\circ\pi(u)-u_n)+\sum_{k=1}^{n-1}\left((g\circ \pi (u)-u_n)D_k g(\pi(u))+\pi_k(u)-u_k\right)D_n \pi_k(u)}{f(u)}\stackrel{\eqref{extremal}}{=}$$
$$-\frac{g\circ \pi(u)-u_n}{f(u)}$$
and, consequently, $|\nabla_u f|=1$ provided that $f(u)>0$, that is $u$ does not belong to the graph of the function $g$, and the partial derivatives of the projection $\pi$ exist at $u\in \mathbb{R}^n$. As a particular e\-xample, consider the graph of the function $g\colon \mathbb{R}\to \mathbb{R}$, $\ g(x):=x^2$. The distance function is given by the formula 
$$f\colon \mathbb{R}^2\to \mathbb{R}, \ f(u):=\inf_{x\in \mathbb{R}}\sqrt{(u_1-x)^2+(u_2-x^2)^2}.$$

\begin{figure}
\centering
\includegraphics[bb=0 0 288 288, scale=0.7]{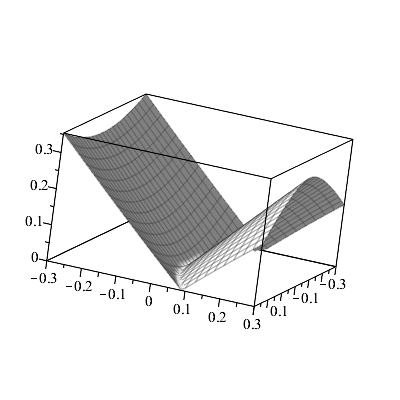}
\caption{The distance function generated by a parabola.}
\end{figure}

\noindent
Since we have a non-negative expression to minimize, it is enough to minimize its square. Differentiating with respect to $x$, the minimizer belonging to $u\in \mathbb{R}^2$ is a real root of the cubic equation
$$u_1-x+2x(u_2-x^2)=0,$$
that is, 
$$x^3+\left(\frac{1}{2}-u_2\right)x-\frac{u_1}{2}=0.$$
If the discriminant
$$D(u):=\frac{u_1^2}{16}+\frac{1}{27}\left(\frac{1}{2}-u_2\right)^3$$
is positive, then we have exactly one real solution
$$\pi(u):=\sqrt[3]{\frac{u_1}{4}+\sqrt{D(u)}}+\sqrt[3]{\frac{u_1}{4}-\sqrt{D(u)}}$$
depending on $u\in \mathbb{R}^2$. Therefore
$$f(u)=\sqrt{(u_1-\pi(u))^2+(u_2-\pi^2(u))^2}$$
provided that $D(u)>0$. The differentiability of the distance function is guaranteed by the conditions $f(u)>0$, that is $u$ does not belong to the graph of the function $g$ and $u_1\neq \pm 4\sqrt{D(u)}$. Figure 1. shows the distance function around the origin. The singularities along the graph of $g$ admit alternating convex and concave parts. 

In what follows, we are interested in functions without singularities. As a main result, we prove that a continuously dif\-ferentiable convex function with gradient of constant norm must be affine. There are also more general results under stronger regula\-ri\-ty conditions: {\emph{if a twice continuously differentiable function on $\mathbb{R}^n$ has gradient of constant norm everywhere, then the function must be affine}} \cite{Letac} (p. 400). In contrast, Theorem \ref{main} reduces the re\-gu\-larity condition by the assumption of convexity. It is a na\-tural question whether the convexity hypothesis can be dropped by using an alternative approach such as 
\begin{verbatim}
https://math.stackexchange.com/questions/867992
\end{verbatim}

The result can also be stated for smooth functions defined on a complete connected Riemannian space provided that the Ricci curvature is non\-ne\-gative \cite{Sakai}. In addition, the existence of affine functions on the base manifold implies that it is isometric to the Riemannian product $f^{-1}(0)\times \mathbb{R}$.  The question naturally arises whether the proof of Theorem \ref{main} can work in the generality of Riemannian spaces.  

\begin{theorem}
\label{main}
Let $f\colon \mathbb{R}^n\to \mathbb{R}$ be a continuously differentiable convex function. If the gradient is of constant norm then the function is an affine mapping. 
\end{theorem}

Following the steps of the proof, we sketch a possible generalization of the result to functions defined on Hilbert spaces.

\section{The proof of Theorem \ref{main}}

As we shall see, the discussion of single-variable functions is relatively simple, but the multivariate case requires more advanced tools. The proof is based on Brouwer's fixed-point theorem \cite{Brouwer}, see also \cite{Hatcher} (Cor. 2.15): \textit{any continuous mapping of a closed Euclidean ball into itself has a fixed point.} The proof also relies on the Cauchy--Bunyakovsky--Schwarz inequality and the first-order characterization of convexity:
\begin{equation}
\label{1st}
f(v)-f(u)\geq \langle \nabla_u f, v-u\rangle, \quad u, v \in \mathbb{R}^n.
\end{equation}

The proof becomes complete by minimizing the distance between lines of gradient directions. 

\subsection{The case of single variable functions} The case of $n=1$ is easy to see because the condition for the gradient reduces to the simple differential equation $f'=\textrm{\ const.}$, where the absolute value is omitted due to the continuity of the derivative. The solutions are of the form $f(u)=c_1 u+c_0$ for some real constants $c_1$, $c_0\in \mathbb{R}$.  

\subsection{The case of multivariable functions} Let $n\geq 2$ and suppose that the gradient is of constant norm. If it is zero then we are done (the function is constant). Otherwise, using a vertical scaling, we can suppose that the gradient is of constant unit norm. Let $r>0$ be a positive real number, $B_r=\{ u\in \mathbb{R}^n \ | \ |u|\leq r \}$ and consider the continuous mapping
$$u\in B_r \mapsto r\nabla_u f$$
of the ball $B_r$ into itself. Using Brouwer's fixed-point theorem we have a one-parameter family $u_r$ of fixed points such that
\begin{equation}
\label{eq:01}
u_r=r\nabla_{u_r}f \quad \Rightarrow \quad u_r/r=\nabla_{u_r} f \quad \Rightarrow \quad \lim_{r\to 0^+} u_r/r=\nabla_{\mathbf{0}}f.
\end{equation}
Taking the one-parameter family of functions $\varphi_r(s):=f(su_r/r)$ we have that
\begin{equation}
\label{eq:02}
\varphi_r'(s)=\langle \nabla_{su_r/r}f, u_r/r\rangle \leq 1
\end{equation}
by the Cauchy-Bunyakovsky-Schwarz inequality. Since $\varphi_r'(r)=1$ and, by the convexity, $\varphi_r'$ is monotone increasing, it follows that
\begin{equation}
\label{eq:03}
\varphi_r'(s)=1 \quad (s\geq r) \quad \Rightarrow \quad \varphi_r(s)=s+f(u_r) \quad (s\geq r).
\end{equation}
Using a continuity argument, the limit $r\to 0^+$ gives that
\begin{equation}
\label{eq:04}
f\left(s \nabla_{\mathbf{0}} f\right)=s+f(\mathbf{0}) \quad (s\geq 0).
\end{equation}
Changing the origin by  a horizontal translation we can reproduce formula \eqref{eq:04} for any $u\in \mathbb{R}^n$, that is
\begin{equation}
\label{eq:05}
\varphi^u(s):=f\left(u+s \nabla_u f\right)=s+f(u) \quad (u\in \mathbb{R}^n, s\geq 0).
\end{equation}
Let us apply formula \eqref{eq:05} to $-v_0$, where $v_0$ is a fixed point of the mapping
$$v\in B_1 \mapsto \nabla_{-v}f \in B_1 \quad \Rightarrow \quad v_0=\nabla_{-v_0}f.$$
We have that 
\begin{gather}
\label{eq:06}
\varphi^{-v_0}(s)=f\left(-v_0+s \nabla_{-v_0} f\right)=s+f(-v_0) \quad (s\geq 0) \quad \Rightarrow \quad\\
f(\mathbf{0})=\varphi^{-v_0}(1)=1+f(-v_0)\notag
\end{gather}
and
\begin{equation}
\label{eq:07}
1=(\varphi^{-v_0})'(1)=\langle \nabla_{\mathbf{0}}f, \nabla_{-v_0} f\rangle \leq 1.
\end{equation}
By the Cauchy-Bunyakovsky-Schwarz inequality,  
$$\nabla_{\mathbf{0}}f=\nabla_{-v_0} f=v_0$$
and we have that 
\begin{equation}
\label{eq:08}
f\left(s\nabla_{\mathbf{0}}f\right)=f\left(-v_0+(s+1) \nabla_{-v_0} f\right)\stackrel{\eqref{eq:06}}{=}(s+1)+f(-v_0)=s+f(\mathbf{0})
\end{equation}
for any $s\geq -1$. Repeating the process it follows that
\begin{equation}
\label{eq:09}
f\left(s\nabla_{\mathbf{0}}f\right)=s+f(\mathbf{0}) \quad (s\in \mathbb{R}).
\end{equation}
Changing the origin by a horizontal translation we can reproduce formula \eqref{eq:09} for any $u\in \mathbb{R}^n$, that is
\begin{equation}
\label{eq:10}
\varphi^u(s):=f\left(u+s \nabla_u f\right)=s+f(u) \quad (u\in \mathbb{R}^n, s\in \mathbb{R}).
\end{equation}
Especially, the gradient is constant along the parameterized line $c_u(s)=u+s\nabla_u f$ because
\begin{equation}
\label{eq:10a}
1=(\varphi^{u})'(s)=\langle \nabla_{u+s\nabla_u f}f, \nabla_{u} f\rangle \leq 1
\end{equation}
and, by the Cauchy-Bunyakovsky-Schwarz inequality, 
\begin{equation}
\label{eq:10b}
\nabla_{u+s\nabla_u f}f=\nabla_{u} f.
\end{equation}
Taking the parameterized lines  $c_u(s)=u+s\nabla_u f$ and $c_v(s)=v+s\nabla_v f$, let $u_0$, $v_0$ be such a pair of points minimizing the distance between them. The difference vector $v_0-u_0$ is obviously orthogonal to the directions $\nabla_u f$ and $\nabla_v f$. Using \eqref{eq:10b}, the first order characterization of the convexity \eqref{1st}  implies that
\begin{gather}
\label{eq:11}
f(v_0)-f(u_0)\geq \langle \nabla_{u_0}f, v_0-u_0\rangle= \langle \nabla_{u}f, v_0-u_0\rangle=0 \quad \Rightarrow \quad \\
f(v_0)\geq f(u_0) \notag.
\end{gather}
 Changing the role of $u_0$ and $v_0$, we have that $f(v_0)=f(u_0)$. Moreover, by equation \eqref{eq:10},
\begin{equation}
\label{eq:12}
f\left(u_0+s \nabla_{u_0} f\right)=s+f(u_0)=s+f(v_0)=f\left(v_0+s \nabla_{v_0} f\right) \quad (s\in \mathbb{R}).
\end{equation}
Therefore, by the first-order characterization \eqref{1st} of convexity,
\begin{gather}
\label{eq:13}
0=f(v_0+s\nabla_{v_0} f)-f(u_0+s\nabla_{u_0}f)\geq \langle \nabla_{u_0}f, v_0-u_0\rangle+ \\
s \langle \nabla_{u_0}f, \nabla_{v_0}f-\nabla_{u_0}f\rangle=s \langle \nabla_{u_0}f, \nabla_{v_0}f-\nabla_{u_0}f\rangle \quad (s\in \mathbb{R}) \notag
\end{gather}
because the gradient is constant along the parameterized lines of the form $c_u(s)=u+s\nabla_u f$ in the sense of equation \eqref{eq:10b}, that is
$$\nabla_{u_0}f=\nabla_{u_0+s\nabla_{u_0}f}f \quad (s\in \mathbb{R})$$
and $\nabla_{u_0}f \ \bot \ v_0-u_0$ because of the choice of $u_0$ and $v_0$. We have 
$$\langle \nabla_{u_0}f, \nabla_{v_0}f-\nabla_{u_0}f\rangle=0 \quad \Rightarrow \quad \langle \nabla_{u_0}f, \nabla_{v_0}f\rangle=1$$
and, by the Cauchy-Bunyakovsky-Schwarz inequality,
$$\nabla_{u_0}f=\nabla_{v_0}f \quad \Rightarrow \quad \nabla_{u}f=\nabla_{v}f.$$
In summary, the gradient is a constant unit vector and
$$f(u)=\langle c_1, u\rangle+c_0 \quad (u\in \mathbb{R}^n),$$
where $c_1\in \mathbb{R}^n$ is a given unit vector (gradient) and $c_0$ is a real constant. 

\begin{corollary}
Let $f\colon \mathbb{R}^n\to \mathbb{R}$ be a continuously differentiable concave function. If the gradient is of constant norm then the function is an affine mapping. 
\end{corollary}

\begin{corollary}
Let $f\colon \mathbb{R}^n\to \mathbb{R}$ be a continuously differentiable convex/concave function and suppose that it is bounded from below/above. If the gradient is of constant norm then the function is constant. 
\end{corollary}

\subsection{The case of Hilbert spaces.} Theorem \ref{main} can also be stated for functions $f\colon H\to \mathbb{R}$ defined on separable Hilbert spaces. Instead of Brouwer's fixed-point theorem we can use the theorem due to F. Browder and J. Minty \cite{Browder}, see also \cite{Renardi} (Thm. 10.49): {\emph{a bounded, continuous, coercive and monotone mapping of a real, separable reflexive Banach space into its continuous dual space is surjective}}. In case of a separable Hilbert space it can be reformulated as follows: {\emph{a bounded, continuous, coercive and monotone mapping of a real, separable Hilbert space into itself is surjective}}. The proof of Theorem \ref{main} works with a slight modification as follows. 

\begin{theorem}
\label{main01}
Let $H$ be a real, separable Hilbert space and consider a continuously differentiable convex function $f\colon H\to \mathbb{R}$. If the gradient is of constant norm, then the function is an affine mapping. 
\end{theorem}

Let $H$ be a real, separable Hilbert space. Since the gradient map of a continuously differentiable convex function is continuous and monotone in the sense that
$$\langle \nabla_u f-\nabla_v f, u-v\rangle \geq 0$$
for any $u, v\in H$, it follows that the perturbed mapping $u\in H \mapsto u+\nabla_u f$ is continuous, monotone and coercive, that is,
$$\lim_{|u|\to \infty} \frac{\langle u+\nabla_u f, u \rangle}{|u|}=\lim_{|u|\to \infty} |u|+ \frac{\langle \nabla_u f-\nabla_{\bf 0}f, u-{\bf 0} \rangle}{|u|}+\frac{\langle \nabla_{\bf 0} f, u \rangle}{|u|}=\infty.$$
Using that the gradient is of constant norm, we also have that the perturbed mapping is bounded, that is, the image of a bounded set is bounded. Therefore the perturbed mapping satisfies the conditions in the Browder-Minty theorem.

\subsection{The sketch of the proof of Theorem \ref{main01}.} In what follows we are going to reproduce the steps in the proof of Theorem \ref{main}. At first, let $r>0$ be a positive real number. Applying the Browder-Minty theorem to each mapping of the form $u\in H\mapsto u+r\nabla_u f$, we can construct a family of points $u_r\in H$ such that $u_r+r\nabla_{u_r} f={\bf 0}$, that is, $u_r=-r\nabla_{u_r}f$. Since $\lim_{r\to 0^+} u_r=0$ and $\nabla f$ is continuous, we have that
\begin{equation}
\label{eq:01H}
-\nabla_{\mathbf{0}}f=\lim_{r\to 0^+} \frac{u_r}{r}.
\end{equation}
Taking the family of functions $\varphi_r(s):=f(su_r/r)$ on $\mathbb{R}$, it follows that
\begin{equation}
\label{eq:02H}
\varphi_r'(s)=\langle \nabla_{su_r/r}f, u_r/r\rangle \geq -1
\end{equation}
by the Cauchy--Bunyakovsky--Schwarz inequality. Since $\varphi_r'(r)=-1$ and, by the convexity, $\varphi_r'$ is monotone increasing, it follows that $\varphi_r'(s)=-1$ for any $s\leq r$, which implies that 
\begin{equation}
\label{eq:03H}
\varphi_r(s)=-s+f(u_r), \quad s\leq r.
\end{equation}
Using a continuity argument, the limit $r\to 0^+$ gives that
\begin{equation}
\label{eq:04H}
f\left(-s \nabla_{\mathbf{0}} f\right)=-s+f(\mathbf{0}), \quad s\leq 0,
\end{equation}
which is equivalent to equation \eqref{eq:04} and we have equation \eqref{eq:05}. 

In the second step we apply the Browder-Minty theorem to the perturbed gradient of the function $u\in H\mapsto f(-u)$ to give an element satisfying $v_0-\nabla_{-v_0}f={\bf 0}$, that is, $v_0=\nabla_{-v_0}f$. Following the proof of Theorem \ref{main} step by step, we have equation \eqref{eq:10}. 
Especially, the gradient is constant along the parameterized line $c_u(s)=u+s\nabla_u f$ (see equation \eqref{eq:10b}). Therefore the last steps of the proof are taken in finite dimensional subspaces.

\end{document}